# A New Approach to Solve the Fuzzy Nonlinear Unconstraint Optimization Problem

Paresh Kumar Panigrahi[1] and Sukanta Nayak[2] and Sudipta Priyadarshini[3]

[1, 2, 3]Department of Mathematics, School of Advanced Sciences, VIT-AP University, Amaravati, Andhra Pradesh

**Abstract**

This paper investigates fuzzy nonlinear system equations using an optimization approach. Here, the inner-outer direct search technique is used with fuzzy coefficients and vectors to quantify the uncertain solution. The fuzzy nonlinear system of equations is converted into an unconstrained fuzzy multivariable optimization problem with preserving the operating constraints. Then developed fuzzy inner-outer direct search is applied. The developed algorithm is demonstrated through an example problem. The obtained uncertain solutions are depicted and compared with the special case. The results are found to be in good agreement when the membership values are unity.

**Keywords**: Fuzzy set, Triangular Fuzzy Number, Unconstrained Minimization Problem, IODS technique.

## 1. Introduction

The system of equations is one of the important mathematical phenomena that occur frequently in engineering and science. Generally, the numerical methods transform the governing differential equations into algebraic equations. Investigation of real-life problems needs numerical methods which provide a system of nonlinear equations [1]. As finding the solution to a nonlinear system of equations with the measure of field variables with an exact method is challenging work, numerical approaches are adopted to convert the nonlinear system equations to algebraic equations [2]. In addition to this, the involvement of uncertainty makes it more difficult to study the system. The role of uncertainties cannot be avoided in the contest better understanding of the system and field variables. Hence, we have considered the fuzzy uncertainties with the nonlinear system of equations. These uncertainties arise as a consequence of experimental error, impreciseness, and insufficient knowledge of the involved parameters ( [3] [4]). The inclusion of uncertainties results in the algebraic equation to a fuzzy system of nonlinear equations. ( [5] [6] [7]).

The notions of fuzzy set theory and its computation, fuzzy number viz. Triangular Fuzzy Number was introduced in [8]. The fuzzy nonlinear system of equations is one of the applications where TFN can be used to compute the fuzzy unknowns of the system. A general model of a $n \times n$ system of equations using TFNs is discussed here. In this context, ( [9] [10]), the researchers used a recursive subdivision to solve a nonlinear system of polynomial equations expressed on Bernstein's basis.

The preceding research demonstrates that direct and/or numerical ways of investigating fuzzy nonlinear equations may be applied. However, the problem comes when one needs to identify the gradient or slopes. So, direct search optimization techniques can be employed for the same. Therefore, the fuzzy nonlinear systems of equations need to be transformed into unconstrained minimization problems that may be addressed using different strategies. Many articles have provided different ideas and techniques to optimize an unconstrained minimization problem in this area. Several of the techniques are described in ( [7] [11]). Moreover, numerous researchers presented various optimization techniques for the same. In this regard, [12] suggested a hybrid technique for solving a system of nonlinear equations that combines a chaotic optimization algorithm with the quasi-Newton method. The basic idea behind this approach is to look for an initial estimate that will reach the convergent areas of the quasi-Newton method. Besides this, the computation cost may rise owing to the quasi-Newton approach, necessitating the development of a new methodology. A recent paper, [13] considered an optimization challenge regarding for the work strategy at the gate of the Barcelona container port. The author then structured the identical issue with the typical traffic flow as a linear programming problem (LPP), given the limitations acquired from the gate data. The simplex method is then used to solve the problem, and the solution is achieved. After determining the best solution, a sensitivity analysis is created. The author investigated two items: feasibility, which indicates that the solution meets all requirements, and optimality, which indicates that the solution is optimal for at least one characteristic model. The basic stability technique is used to modify the Two-Step Method (TSM). Also, the author provided a technique for solving fully fuzzy linear programming (FFLP) problems using the ranking function. In [14] authors devised an approach for solving convex quadratic programming problems with limited variables. While in ( [15], [16], [17]) researchers proposed a technique for calculating the response interval of structures caused by huge uncertainty( [18] [19]).

---

[1] Email: jrpareshkp@gmail.com (Paresh Kumar Panigrahi)
[2] Email: sukantgacr@gmail.com (Sukanta Nayak)
[3] Email: sudipta.grc@gmail.com (Sudipta Priyadarshini)

As per the above literature review, it is found that there is a need for an alternate approach that can be used to solve the fuzzy nonlinear system of equations with an optimization approach. In this respect, we developed an inner-outer Direct search technique (IODS) direct search optimization algorithm for fuzzy nonlinear systems of equations. The suggested optimization approach transforms a fuzzy system nonlinear of equations into an unconstrained multivariable optimization problem with triangular fuzzy uncertainty in three different steps. Then, the inner-outer direct search strategy is used to examine the fuzzy unconstrained multivariable optimization problem. Finally, the discovered solutions are integrated, and the regularity principle is adopted to get the approximated solutions.

2. **Fuzzy System of Nonlinear Equations**

This section presents a nonlinear system of equations with fuzzy environment. Here, Triangular Fuzzy Numbers (TFNs) are considered for uncertainness and then the fuzzy system of nonlinear equations is depicted. For the sake of completeness first the definition of TFN is discussed. Then the same is incorporated into system of nonlinear equations.

A fuzzy number of $A = [x_L, x_N, x_R]$ is said to be TFN if the membership values can be defined as

$$\mu_A(x) = \begin{cases} f_L, & x_L \leq x \leq x_N \\ f_R, & x_N \leq x \leq x_R \\ 0, & otherwise \end{cases} \quad (1)$$

where, $f_L$ is defined as the left monotonically increasing function and $f_R$ is defined as the right monotonically decreasing function. The functions are given by $f_L = \frac{x - x_L}{x_N - x_L}$ and $f_R = \frac{x_R - x}{x_R - x_N}$.

For computational purpose, an arbitrary TFN $A = [x_L, x_N, x_R]$ can be transformed to multi-variable form using the following two steps.

1. TFN to interval form,
2. Interval to crisp form.

Using $\alpha$-cuts the TFN can be written as

$$\tilde{A} = [x_L, x_N, x_R] \approx [\xi_L, \xi_R], \quad (2)$$

where,

$$\xi_L = x_L + \alpha(x_N - x_L) \text{ and } \xi_R = x_R + \alpha(x_N - x_R); \alpha \in [0,1].$$

The crisp representation of TFN $\tilde{A}$ is

$$\xi_L - \beta(\xi_L - \xi_R), \beta \in [0,1]. \tag{3}$$

The above representation may be used in the following system of equations to compute fuzzy system of equations.

Consider a system of equations in matrix form

$$\tilde{M}x = \tilde{b} \tag{4}$$

where, $\tilde{M}$ is a coefficient matrix, $\tilde{b}$ is the right-side vector, and $x$ is the unknown vector that need to quantify.

The compact form of the same is defined as

$$\tilde{M} = [\tilde{a}_{ij}], i,j = 1,2,\ldots n, \tilde{b} = [\tilde{b}_1, \tilde{b}_2, \ldots, \tilde{b}_n]^T, \text{ and } \tilde{x} = [\tilde{x}_1, \tilde{x}_2, \ldots \tilde{x}_n]^T \tag{5}$$

Here, the entries $\tilde{a}_{ij}$ and $\tilde{b}_j$ are the TFNs.

Consider the system of the equation having $n$ variables

$$\tilde{M}_{n \times n} \tilde{X}_{n \times 1} = \tilde{b}_{n \times 1} \tag{6}$$

Where

$$\tilde{M}_{n \times n} = \begin{bmatrix} \tilde{a}_{11} & \tilde{a}_{12} & \cdots & \tilde{a}_{1n} \\ \tilde{a}_{21} & \tilde{a}_{22} & \cdots & \tilde{a}_{2n} \\ \vdots & \vdots & \ddots & \vdots \\ \tilde{a}_{n1} & \tilde{a}_{n2} & \cdots & \tilde{a}_{nn} \end{bmatrix}, X_{n \times n} = \begin{bmatrix} x_1 \\ x_2 \\ \vdots \\ x_n \end{bmatrix}, \text{ and } \tilde{b}_{n \times n} = \begin{bmatrix} \tilde{b}_1 \\ \tilde{b}_2 \\ \vdots \\ \tilde{b}_n \end{bmatrix}. \tag{7}$$

The expanded form of TFNs $\tilde{a}_{ij}$ and $\tilde{b}_j$ are

$$\tilde{a}_{ij} = [a_{L_{ij}}, a_{N_{ij}}, a_{R_{ij}}] \text{ and } \tilde{b}_j = [b_{L_j}, b_{N_j}, b_{R_j}]; i,j = 1,2,\cdots,n.$$

The above formulation is used in Inner Outer Direct Search (IODS) Optimization Technique and extended the same for fuzzy environment. As such, the next section explains the IODS optimization method for solving fuzzy nonlinear systems of equations.

3. **Inner Outer Direct Search (IODS) Optimization Technique in Fuzzy Environment**

To solve a fuzzy system of nonlinear equations, the IODS optimization technique is applied here. First, the fuzzy nonlinear system of equations is converted to a fuzzy unconstrained optimization problem and then using the extended IODS fuzzy nonlinear system of equations is solved.

The converted fuzzy unconstrained optimization problem produces a minimization problem that minimizes fuzzy objective function. The IODS method is executed though exploratory and pattern search. The search is carried out by using the idea of inner and outer calculations, as well as for the converted multi-variable function where the variable ranges from zero to unity.

For understanding, we have considered a system of two equations, where the right-hand side values are TFN.

$$f_k(x_1, x_2) = [b_{L_k}, b_{N_k}, b_{R_k}] \tag{8}$$

The eq. (6) can be transformed into the interval by using the following parameter form

$$f_k(x_1, x_2) = [\xi_{L_k}, \xi_{R_k}] \tag{9}$$

where $0 \leq \alpha_K \leq 1$

$$\xi_{L_k} = b_{L_k} + \alpha_k(b_{N_k} - b_{L_k}), \xi_{R_k} = b_{L_k} + \alpha_k(b_{N_k} - b_{L_k}), (k = 1,2 \ldots, n)$$

The value of $0 \leq \alpha < 1$ of eq. (8), the left-hand side outer computation can be described in this matrix from

$$\begin{bmatrix} a_{L_{11}} & a_{L_{12}} & \cdots & a_{L_{1n}} \\ a_{L_{21}} & a_{L_{22}} & \cdots & a_{L_{2n}} \\ \vdots & \vdots & \ddots & \vdots \\ a_{L_{n1}} & a_{L_{n2}} & \cdots & a_{L_{nn}} \end{bmatrix} \begin{bmatrix} x_1 \\ x_2 \\ \vdots \\ x_n \end{bmatrix} = \begin{bmatrix} b_{L_1} \\ b_{L_2} \\ \vdots \\ b_{L_n} \end{bmatrix} \quad (10)$$

The value of $0 \leq \alpha < 1$ of eq. (8), the Right-hand side outer computation can be described in this matrix from

$$\begin{bmatrix} a_{R_{11}} & a_{R_{12}} & \cdots & a_{R_{1n}} \\ a_{R_{21}} & a_{R_{22}} & \cdots & a_{R_{2n}} \\ \vdots & \vdots & \ddots & \vdots \\ a_{R_{n1}} & a_{R_{n2}} & \cdots & a_{R_{nn}} \end{bmatrix} \begin{bmatrix} x_1 \\ x_2 \\ \vdots \\ x_n \end{bmatrix} = \begin{bmatrix} b_{R_1} \\ b_{R_2} \\ \vdots \\ b_{R_n} \end{bmatrix} \quad (11)$$

Similarly, using the value of $\alpha = 1$ eq. (8) the inner computation can be represented in this matrix from

$$\begin{bmatrix} a_{N_{11}} & a_{N_{12}} & \cdots & a_{N_{1n}} \\ a_{N_{21}} & a_{N_{22}} & \cdots & a_{N_{2n}} \\ \vdots & \vdots & \ddots & \vdots \\ a_{N_{n1}} & a_{N_{n2}} & \cdots & a_{N_{nn}} \end{bmatrix} \begin{bmatrix} x_1 \\ x_2 \\ \vdots \\ x_n \end{bmatrix} = \begin{bmatrix} b_{N_1} \\ b_{N_2} \\ \vdots \\ b_{N_n} \end{bmatrix} \quad (12)$$

All the left outer, right outer, and inner computation are performed in the given algorithm. The above procedure can be shown in the Inner Outer Search Algorithm projected below.

**Fuzzy IODS algorithm**

Step 1: Construct the left and right outer fuzzy system of the nonlinear equation.

Step 2: Convert the inner and outer fuzzy system of the nonlinear equation into an unconstrained minimization problem.

Step 3: Perform exploratory move

　　Step 3.1:

　　Calculate $F = F(x), F^+ = F(x_i + \Delta_i)$ and $F^- = F(x_i - \Delta_i)$

　　Step 3.2:

　　Find $F_{min} = \min(F, F^+, F^-)$ Set $x$ resembles to $F_{min}$.

　　Step 3.3:

　　Is $i = N$ (Number of independent variables)? If no, set $i = i + 1$ and

go to Step 3.1 Else $x$ is the result and go to Step 3.4

　　Step 3.4:

　　If $x \neq x^c$, success;

　　Else failure.

Step 4: After that, Perform the following pattern search move

　　　Step 4.1: Choose a starting point $x^{(0)}$, variable increments $\Delta_i$ $(i = 1,2, \cdots, N)$, a step reduction factor $a > 1$, and a termination parameter, $\varepsilon$. Set $k = 0$.

Step 4.2: Perform an exploratory move using $x^{(k)}$ as the base point. consider $x$ is the result of the exploratory move. If the exploratory move is a successful, set $x^{(k+1)} = x$ and proceed to Step4.4.

Else go to Step 4.3.

Step 4.3: Is $\|\Delta\| < \epsilon$ ? If yes, Terminate

Else set $\Delta_i = \frac{\Delta_i}{a}$ for $i = 1, 2, \cdots, N$ and go to Step 4.2.

Step 4.4: Set $k = k + 1$ and perform the pattern move

$$x_p^{(k+1)} = x^{(k)} + (x^{(k)} - x^{(k-1)}).$$

Step 4.5: Perform another exploratory move with $x_p^{(k+1)}$ as the best point. Let the result be $x_p^{(k+1)}$.

Step 4.6: Is $f(x^{(k+1)}) < f(x^{(k)})$? If yes, go to Step 4.4. Else go to Step 4.3.

Step 5: Assemble all the solutions and sort the minimum and maximum individual components of the solution vector.

Using this algorithm, one need to solve the fuzzy unconstrained optimization problem.

The solution is taken as three values of left side outer, inner and right side outer as $x_L^K, x_N^K$ and $x_R^K$, where $k$ is the number of variables ($k = 1, 2, 3, \ldots$). Then the final three solutions need to be arranged to get the desired solution.

4. **Test Problem**

Here, we have considered the following example to demonstrate the proposed algorithm.

Consider a fuzzy system of nonlinear equations

$$\begin{aligned} x_1^2 + x_2 &= [2, 5, 8] \\ x_1^2 + x_2^2 &= [3, 6, 9] \end{aligned} \tag{13}$$

We take the initial approximation $x^0 = (1,1)^T$, step size $\Delta = (0.5, 0.5)^T$, and $\epsilon = 0.001$.

The first step is to construct the outer system and then convert Eq (13) into a fuzzy unconstrained minimalization problem. The transformed fuzzy unconstrained minimization problem statement is

$$F(x) = (x_1^2 + x_2 - [2, 5, 8])^2 + (x_1^2 + x_2^2 - [3, 6, 9])^2 \tag{14}$$

Now applying the IODS algorithm, the desired solution is obtained after 19 iterations using MATLAB. The solution vector after 19 iterations is $(\tilde{x}_1, \tilde{x}_2)$, where, $\tilde{x}_1 = [0.6938, 1.7115, 2.3274]$ and $\tilde{x}_2 = [1.5186, 2.0709, 2.5831]$. The solution components $\tilde{x}_1$ and $\tilde{x}_2$ of the solution vector are shown in Fig. 1 and 2 respectively.

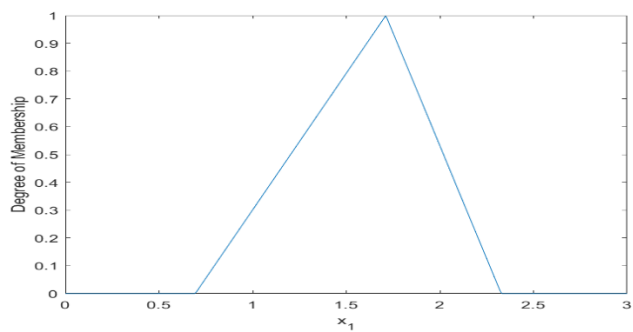

Fig. 1 The uncertain of solution of $x_1$ TFN

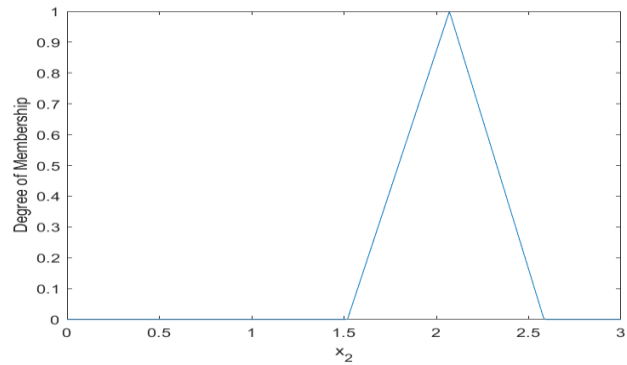

Fig. 1 The uncertain of solution of $x_2$ TFN

### 5. Conclusion

In this paper, the fuzzy IODS optimization approach is proposed to solve fuzzy system of nonlinear equations. The fuzzy IODS method is very simple to use, but before applying it to the fuzzy system of the nonlinear equation, it must be converted to the fuzzy unconstrained optimization problem. Then, triangular fuzzy numbers are used to substitute converted single variable functions. To find the minimum of the fuzzy objective function, the approach employs exploratory and pattern search motions, as well as the idea of inner and outer computation. Finally, the three solution combinations are found and ordered using the regularity principle to achieve the required solution. The computational advantage and simple methodology make it convenient to apply.